\documentclass[12pt]{amsart}

\usepackage{amssymb,amscd,amsmath,amsthm}

\allowdisplaybreaks

\theoremstyle{plain}
\newtheorem{thm}{Theorem}

\newtheorem{lem}[thm]{Lemma}

\newtheorem{cor}[thm]{Corollary}
\newtheorem{prop}[thm]{Proposition}

\theoremstyle{definition}

\theoremstyle{remark}

\newtheorem*{remark*}{Remark}

\newcommand{\cF}{{\mathcal{F}}}

\newcommand{\cJ}{{\mathcal{J}}}

\newcommand{\cV}{{\mathcal{V}}}
\newcommand{\cW}{{\mathcal{W}}}

        \newcommand{\field}[1]{{\mathbb{#1}}}
        \newcommand{\N}{\field{N}}
        \newcommand{\Z}{\field{Z}}
        
        \newcommand{\R}{\field{R}}
        \newcommand{\C}{\field{C}}

\newcommand{\supp}{\operatorname{supp}}

\newcommand{\Dom}{\mbox{\rm Dom}}

\newcommand{\id}{\mbox{\rm id}}
\newcommand{\Tr}{\mbox{\rm Tr}}

\begin{document}
\bibliographystyle{plain}

\title[Periodic Schr\"odinger operators with magnetic wells]
{Semiclassical asymptotics and gaps in the spectra of periodic
Schr\"odinger operators with magnetic wells}

\author{Bernard Helffer}
\address{D\'epartement de Math\'ematiques, B\^atiment 225, F91405 Orsay C\'edex}
\email{Bernard.Helffer@math.u-psud.fr}

\author{Yuri A. Kordyukov}
\address{Institute of Mathematics, Russian Academy of Sciences, 112 Chernyshevsky str., 450077
Ufa, Russia} \email{yuri@imat.rb.ru}
\thanks{B.H. acknowledges support  from the SPECT programme
 of the ESF and from  the European Research Network
`Postdoctoral Training Program in Mathematical Analysis of Large
Quantum Systems' with contract number HPRN-CT-2002-00277. Y.K. acknowledges support from Russian Foundation of Basic Research
(grant no. 04-01-00190).}

\begin{abstract}
We show that, under some very weak assumption of effective variation
for  the
magnetic field, a periodic Schr\"odinger operator with magnetic
wells on a noncompact Riemannian manifold $M$ such that $H^1(M,
\R)=0$ equipped with a properly disconnected, cocompact action of
a finitely generated, discrete group of isometries has an
arbitrarily large number of spectral gaps in the semi-classical
limit.
\end{abstract}
\maketitle

\section{Introduction}
Let $ M$ be a noncompact oriented manifold of dimension $n\geq 2$
equipped with a properly disconnected action of a finitely
generated, discrete group $\Gamma$ such that $M/\Gamma$ is
compact. Suppose that $H^1(M, \R) = 0$, i.e. any closed $1$-form
on $M$ is exact. Let $g$ be a $\Gamma$-invariant Riemannian metric
and $\bf B$ a real-valued $\Gamma$-invariant closed 2-form on $M$.
Assume that $\bf B$ is exact and choose a real-valued 1-form $\bf
A$ on $M$ such that $d{\bf A} = \bf B$.

Consider a Schr\"odinger operator with magnetic potential $\bf A$,
\[
H^h = (ih\,d+{\bf A})^* (ih\,d+{\bf A}),
\]
as a self-adjoint operator in the Hilbert space $L^2(M)$. Here
$h>0$ is a semiclassical parameter, which is assumed to be small.

For any $x\in M$ denote by $B(x)$ the anti-symmetric linear
operator on the tangent space $T_x{ M}$ associated with the 2-form
$\bf B$:
\[
g_x(B(x)u,v)={\bf B}_x(u,v),\quad u,v\in T_x{ M}.
\]
Recall that the intensity of the magnetic field is defined as
\[
{\Tr}^+ (B(x))=\sum_{\substack{\lambda_j(x)>0\\ i\lambda_j(x)\in
\sigma(B(x)) }}\lambda_j(x)=\frac{1}{2}\Tr([B^*(x)\cdot
B(x)]^{1/2}).
\]
Let
\[
b_0=\min \{{\Tr}^+ (B(x))\, :\, x\in { M}\}.
\]
We will always assume that there exist a (connected) fundamental
domain $\cF$ and $\epsilon_0>0$ such that
\begin{equation}\label{e:tr1}
  {\Tr}^+ (B(x)) \geq b_0+\epsilon_0, \quad x\in \partial\cF.
\end{equation}
For any $\epsilon_1 \leq \epsilon_0$, let
\[
U_{\epsilon_1} = \{x\in \cF\,:\, {\Tr}^+ (B(x)) < b_0+
\epsilon_1\}.
\]
Thus $U_{\epsilon_1}$ is an open subset of $\cF$ such that
$U_{\epsilon_1}\cap \partial\cF=\emptyset$ and, for $\epsilon_1 <
\epsilon_0$, $\overline{U_{\epsilon_1}}$ is compact and included
in the interior of $\cF$. Any connected component of
$U_{\epsilon_1}$ with $\epsilon_1 < \epsilon_0$ can be understood
as a magnetic well (attached to the effective potential $h\cdot
{\Tr}^+ (B(x))$).

Consider the set $U^+_{\epsilon_0}$, which consists of all $x\in
U_{\epsilon_0}$ such that the rank of $B(x)$ is locally constant
at $x$, that is, constant in an open neighborhood of $x$. Let us
assume that
\begin{equation}\label{e:trB}
  {\Tr}^+ B\ \text{is not locally constant on}\ U^+_{\epsilon_0}.
\end{equation}
The assumption (\ref{e:trB}) holds for any $B$, satisfying the
assumption~(\ref{e:tr1}), if the dimension $n$ equals $2$ or $3$.

If $T$ is a self-adjoint operator, $\sigma(T)$ denotes its
spectrum. By a gap in the spectrum of $T$ we will mean an interval
$(a,b)$ such that
\[
(a,b)\cap \sigma(T) = \emptyset\,.
\]

\begin{thm}\label{t:main}
Under the assumptions (\ref{e:tr1}) and (\ref{e:trB}), there
exists, for any natural $N$, $h_0>0$ such that, for any $h\in
(0,h_0]$, the spectrum of $H^h$ contained in $[0,
h(b_0+\epsilon_0)]$ has at least $N$ gaps.
\end{thm}

The proof of Theorem~\ref{t:main} is based on the study of the
tunneling effect for the operator $H^h$. First, we prove that the
spectrum of $H^h$ is localized inside an exponentially small
neighborhood of the spectrum of its Dirichlet realization $H^h_D$
in $D=\overline{U_{\varepsilon_1}}$ for $\varepsilon_1 <
\varepsilon_0$ (a multi-well problem). For this, we follow the
approach to the study of the tunneling effect in multi-well
problems developed by Helffer and Sj\"ostrand for Schr\"odinger
operators with electric potentials (see for instance
\cite{HSI,HSII}) and extended to magnetic Schr\"odinger operators
in \cite{HS87,HM}. Since $H^h$ is not with compact resolvent, we
work not with individual eigenfunctions as in \cite{HSI}, but with
resolvents, using the strategy developed in
\cite{HSII,HS88,Di-Sj,Carlsson} for the case of electric potential
and in \cite{Frank} for the case of magnetic field. The idea is to
construct an approximate resolvent $R^h(z)$ of the operator $H^h$
for any $z$, which is not exponentially close to the spectrum of
$H^h_D$, starting from the resolvent of $H^h_D$ and the resolvent
of the Dirichlet realization of $H^h$ in the complement to the
wells. The proof of the fact that the error of the approximation
is exponentially small is based on Agmon-type weighted estimates
(cf. \cite{Ag} and their semi-classical versions in \cite{HSI} for
the case of Schr\"odinger operators and \cite{HM} for the case of
magnetic Schr\"odinger operators). A very related result was
proved by Nakamura in \cite{Nakamura95}.

Thus the proof is reduced to the study of the discrete spectrum of
the operator $H^h_D$ in the interval $[0, h(b_0+\epsilon_0)]$.
Actually, it remains to show that there is, as $h\to 0$, an
arbitrarily large number of gaps in the spectrum of $H^h_D$ of
size $> h^M$ with some constant $M>0$. For this, we use a weak
polynomial upper bound on the number of eigenvalues of $H^h_D$
contained in $[0, h(b_0+\epsilon_0)]$ and the construction of
quasimodes of the operator $H^h$ given in the proof of
\cite[Theorem 2.2]{HM}. We need the fact that the number of
quasimodes, which we can construct, is sufficiently large, that
follows from a slight modification of \cite[Theorem 2.2]{HM} (see
Proposition~\ref{t:modes} below) and makes an essential use of the
assumption~(\ref{e:trB}).

It seems that the periodicity assumption in Theorem~\ref{t:main}
is not important, and the theorem could be probably extended to
the case when we only assume the existence of an infinite number
of identical magnetic wells of the form $U_{\epsilon_0}$ separated
by regions where the estimate~(\ref{e:tr1}) holds. To show such a
result, one should use the strategy developed in \cite{Carlsson}
in the case of the strong electric field and in
\cite{Cornean-Nenciu} for the tight binding model, but this will
not be detailed in this article.

Let us mention some related results on gaps in the spectrum of
the operator $H^h$.

An asymptotic description of the spectrum of the
two-di\-men\-si\-o\-nal magnetic Schr\"odinger operator with a
periodic potential in a strong magnetic field can be given, using
averaging methods or effective hamiltonians together with
semiclassical analysis (see, for instance,
\cite{BDP,HS88,HSLNP345} and references therein). This allows to
give, at least, heuristically, a more precise asymptotic picture
of spectral bands and gaps for these operators. However, it
should be noted that, in these papers, the magnetic field is,
usually, supposed to be uniform and spectral gaps are created by
the electric potential, whereas in our case the electric field
vanishes and spectral gaps are created by a periodic array of
magnetic barriers.

In \cite{HempelHerbst95}, Hempel and Herbst studied the strong
magnetic field limit ($\lambda\to\infty$) for the periodic
Schr\"odinger operator in $\R^n$:
\[
H_{\lambda {\bf A},0}=(D-\lambda {\bf A})^2, \quad
D_j=\frac{1}{i}\frac{\partial}{\partial x_j},
\]
where ${\bf B} = d{\bf A}$ is a $\Z^n$-periodic $2$-form. Let
$S=\{x\in\R^n: {\bf B}(x)=0\}$ and $S_{\bf A}=\{x\in\R^n: {\bf
A}(x)=0\}$. Assume that the set $S\setminus S_{{\bf A}}$ has
measure zero, the interior of $S$ is non-empty and $S$ can be
represented as $S=\cup_{j\in {\Z^n}} S_j$ (up to a set of measure
zero) where the $S_j$ are pairwise disjoint compact sets with
$S_j=S_0+j$. It is shown that, as $\lambda\to\infty$, $H_{\lambda
{\bf A},0}$ converges in norm resolvent sense to the Dirichlet
Laplacian $-\Delta_S$ on the closed set $S$. Therefore, as
$\lambda\to\infty$, the spectrum of $H_{\lambda {\bf A},0}$
concentrates around the eigenvalues of $-\Delta_S$ and gaps opens
up in the spectrum of $H_{\lambda {\bf A},0}$. For the operator
$H^h=h^2H_{{h^{-1}\bf A},0}$ this means that for any natural $N$
there exist $C>0$ and $h_0>0$ such that the part of the spectrum
of $H^h$ contained in the interval $[0,Ch^2]$ has at least $N$
spectral gaps for any $h\in (0,h_0)$. The rate of approach of the
resolvent $(H_{\lambda {\bf A},0}+1)^{-1}$ to a limit was studied
by Nakamura in \cite{Nakamura95}.

The case when the set $S\setminus S_{{\bf A}}$ has nonzero measure
was studied by Herbst and Nakamura in \cite{HerbstNakamura}. They
showed that in many situations of interest where this condition
holds the equivalence class of $H_{\lambda {\bf A},0}$ approaches
a periodic or almost-periodic orbit in the space of such classes
as $\lambda\to\infty$, and, therefore, the spectrum of $H_{\lambda
{\bf A},0}$ approaches a periodic or almost-periodic orbit in the
space of subsets of $[0,\infty)$.

In \cite{Ko04}, the author investigated the case when the bottom
$S$ of magnetic wells has measure zero and the magnetic field has
regular behavior near $S$. More precisely, assume that there
exists at least one zero of $B$, and, for some integer $k>0$, if
$B(x_0)=0$, then there exists a positive constant $C$ such that,
for all $x$ in some neighborhood of $x_0$,
\begin{equation}\label{e:B}
C^{-1}|x- x_0|^k\leq {\Tr}^+ (B(x))  \leq C |x- x_0|^k\,.
\end{equation}
It is shown in \cite{Ko04} that, under these assumptions, there
exists an increasing sequence $\{\lambda_m, m\in\N \}$, satisfying
$\lambda_m\to\infty$ as $m\to\infty$, such that, for any $a$ and
$b$, satisfying $\lambda_m<a<b < \lambda_{m+1}$ with some $m$,
\[
[ah^{\frac{2k+2}{k+2}}, bh^{\frac{2k+2}{k+2}}]\cap
\sigma(H^h)=\emptyset\,,
\]
for any $h>0$ small enough. In particular, this implies that, for
any natural $N$, there exist $C>0$ and $h_0>0$ such that the part
of the spectrum of $H^h$ contained in the interval
$[0,Ch^{\frac{2k+2}{k+2}}]$ has at least $N$ gaps for any $h\in
(0,h_0)$.

The results of this paper can be considered as a complement of the
results of \cite{Ko04} and, in some sense, correspond to the case
when the condition (\ref{e:B}) holds with $k=0$ (whereas the
results of \cite{HempelHerbst95} are related with the case when
the condition (\ref{e:B}) holds with arbitrarily large $k$). From
the other side, it should be noted that here we state only the
existence of an arbitrarily large number of spectral gaps in the
semi-classical limit and don't know any results on the spectral
concentration in this case.

The authors are grateful to the referee for useful remarks.

\section{Proof of the main theorem}
For any domain $W$ in $ M$, denote by $H^h_W$ the operator $H^h$
in $\overline{W}$ with the Dirichlet boundary conditions. The
operator $H^h_W$ is generated by the quadratic form
\[
q^h_W [u] : = \int_W |(ih\,d+{\bf A})u|^2\,dx
\]
with the domain
\[
\Dom (q^h_W) = \{ u\in L^2(W) : (ih\,d+{\bf A})u \in
L^2\Omega^1(W), u\left|_{\partial W}\right.=0 \},
\]
where $L^2\Omega^1(W)$ denotes the Hilbert $L^2$ space of
differential $1$-forms on $W$, $dx$ is the Riemannian volume form
on $ M$.

Let us assume that (\ref{e:tr1}) and (\ref{e:trB}) are satisfied.
By (\ref{e:trB}), there exists a connected open set $\Omega\subset
U_{\epsilon_0}$ such that the rank of $B(x)$ is constant on
$\Omega$ and ${\Tr}^+ B(\overline\Omega) =[\alpha,\beta], \alpha
<\beta$. Without loss of generality, we can assume that
$\Omega\subset U_{\epsilon_1}$ for some $\epsilon_1<\epsilon_0$
and, therefore, $[\alpha,\beta]\subset [0, b_0+\epsilon_1]$.

For a fixed $\epsilon_2$ such that $\epsilon_1 < \epsilon_2 <
\epsilon_0$, consider the operator $H^h_D$ associated with
$D=\overline{U_{\epsilon_2}}$. The operator $H^h_D$ has discrete
spectrum. Denote by $\lambda_1^h< \lambda_2^h< \ldots <
\lambda^h_{N(h)}$ the eigenvalues of $H^h_D$ contained in the
interval $[h\alpha,h\beta]$. It follows from rough estimates for
the eigenvalue counting function of $H_D^h$ (cf. for instance
\cite[Lemma 4.2]{HM}) that there exist $C$ and $h_0$ such that
\begin{equation}\label{l:Nh}
  N(h) \leq Ch^{-n}, \quad \forall h\in (0,h_0]\;.
\end{equation}

\begin{thm}\label{t:D}
Under the assumption (\ref{e:tr1}), for any $\epsilon_1 <
\epsilon_2< \epsilon_0$, there exist $C, c, h_0>0$ such that for
any $h\in (0,h_0]$
\[
\sigma(H^h)\cap [0, h(b_0+\epsilon_1)] \subset \{\lambda\in [0,
h(b_0+\epsilon_1)] : {\rm dist}(\lambda, \sigma(H^h_D))<
Ce^{-c/\sqrt{h}}\}.
\]
\end{thm}

The proof of Theorem~\ref{t:D} will be given in
Section~\ref{s:weight}. A slightly weaker version of this theorem
(which involves the largest absolute value of the eigenvalues of
$B(x)$ instead of ${\Tr}^+ (B(x))$) was proved in
\cite{Nakamura95}.

By Theorem~\ref{t:D}, $\sigma(H^h)\cap [h\alpha,h\beta]$ is
contained in exponentially small neighborhoods of $\lambda^h_j,
j=1,2,\cdots, N(h)$: there exist $C, c, h_0>0$ such that for any
$h\in (0,h_0]$
\begin{equation}\label{e:exp}
\sigma(H^h)\cap [h\alpha,h\beta]\subset \bigcup_{j=1}^{N(h)}
[\lambda^h_j-Ce^{-c/\sqrt{h}}, \lambda^h_j+Ce^{-c/\sqrt{h}}].
\end{equation}
It follows from (\ref{e:exp}) that for any $j$ such that
$\lambda^h_{j+1}- \lambda^h_j \geq h^M$ with some $M>0$ the
interval $(\lambda^h_j+Ce^{-c/\sqrt{h}},
\lambda^h_{j+1}-Ce^{-c/\sqrt{h}})$ is a gap in the spectrum of
$H^h$ if $h$ is small enough. Therefore, the proof of
Theorem~\ref{t:main} is completed by the following fact.

\begin{prop}
There exists a constant $M>0$ such that the number of
$j\in\{1,2,\cdots,N(h)-1\}$ with $\lambda^h_{j+1}- \lambda^h_j
\geq h^M$ tends to infinity as $h\to 0$.
\end{prop}

\begin{proof}
First, observe that there exists a constant $C_1>0$ such that, for
any $j=1,2,\cdots,N(h)-1$, we have
\begin{equation}\label{e:11}
\lambda^h_{j+1}- \lambda^h_j \leq C_1h^{4/3},
\end{equation}
and also
\begin{equation}\label{e:22}
\lambda^h_{1}- h\alpha \leq C_1h^{4/3}, \quad h\beta-
\lambda^h_{N(h)} \leq C_1h^{4/3}.
\end{equation}
To see this, we will use the following proposition, which is a
slight modification of \cite[Theorem 2.2]{HM}.

\begin{prop}\label{t:modes}
Assume that the rank of $B$ is constant in a connected open subset
$\Omega$. For any compact subset $K$ of $\Omega$, there exists
$C>0$ such that, for any $\mu$ in ${\Tr}^+ B(K)$ and for any $h\in
(0,1]$
\[
(-h^{4/3}C+h\mu, h\mu + h^{4/3}C)\cap \sigma(H_D^h) \neq \emptyset\,
.
\]
\end{prop}

\begin{proof}
We will follow the proof of \cite[Theorem 2.2]{HM}. Denote by $2d$
the rank of $B(y)$, $y\in \Omega$. By assumption, $d$ is
independent of $y$. For any $y\in \Omega$, there exists an
orthonormal base $e_1(y), e_2(y),\cdots , e_n(y)$ in $T_y { M}$
such that
\begin{align*}
B(y)e_{2j-1}(y)=\mu_j(y)e_{2j}(y), \quad & j=1,2,\cdots, d\,,\\
B(y)e_{2j}(y)=-\mu_j(y)e_{2j-1}(y), \quad & j=1,2,\cdots, d\,,\\
B(y)e_{2d+k}(y)=0, \quad & k=1,2,\cdots, n-2d\,.
\end{align*}
Moreover, for any $j$, $\mu_j(y)$ depends continuously on $y\in
\Omega$, and one can choose the orthonormal base $e_1(y),
e_2(y),\cdots , e_n(y)$, depending continuously on $y$. Let
$\phi_y: \cV_y \to \R^n$ be a local coordinate chart given by the
normal coordinates of the metric $g$ associated with the
orthonormal base $e_1(y), e_2(y),\cdots , e_n(y)$. Without loss of
generality, we can assume that $\phi_y$ is defined in a
neighborhood $\cV_y\subset \Omega$ of $y$ and $\phi_y(\cV_y)$ is a
fixed ball $B$ in $\R^n$ centered at the origin. Moreover, the
family $\{\phi_y^{-1} : y\in \Omega\}$ yields a continuous family
of smooth maps from $B$ to $ M$. In the coordinates $\phi_y$,
$g_{y}$ becomes the standard Euclidean metric on $\R^n$ and
\[
{\bf B}(y)=\sum_{j=1}^d\mu_j(y)\,dx_{2j-1}\wedge dx_{2j}.
\]

Now one can proceed as in the proof of \cite[Theorem 2.2]{HM} and
construct a continuous family $u^h_y\in C^\infty_c(\cV_y)\subset
L^2(D), y\in \Omega,$ such that
\[
\|(H^h-h{\Tr}^+ (B(y)))u^h_y\|_{L^2(D)}\leq
Ch^{4/3}\|u^h_y\|_{L^2(D)}, \quad y\in K,
\]
where $C$ is independent of $y\in K$ by continuity in $y$, that
immediately concludes the proof.
\end{proof}

By Proposition~\ref{t:modes}, the operator $H^h_D$ cannot have
spectral gaps of size, larger than $C_1h^{4/3}$ with some $C_1>0$,
in the interval $[h\alpha, h\beta]$, that immediately implies the
estimates (\ref{e:11}) and (\ref{e:22}).

Now assume from the contrary that for any real $M$ the cardinality
of the set
\[
{\cJ}^h_{M}=\{j\in\{1,2,\cdots,N(h)-1\} : \lambda^h_{j+1}-
\lambda^h_j \geq h^M\}
\]
is bounded as $h\to 0$:
\begin{equation}\label{e:K}
\sharp {\cJ}^h_{M}\leq K, \quad h\in (0,1],
\end{equation}
where $K$ is independent of $h$. Then, using (\ref{l:Nh}),
(\ref{e:11}), (\ref{e:22}) and (\ref{e:K}), we get, for all
sufficiently small $h>0$,
\begin{multline*}
h(\beta-\alpha)=(h\beta-\lambda^h_N)+\sum_{j=1}^{N-1}(\lambda^h_{j+1}-
\lambda^h_j)+ (\lambda^h_1-h\alpha)\\
=(h\beta-\lambda^h_N)+\sum_{j\in {\cJ}^h_{M}}(\lambda^h_{j+1}-
\lambda^h_j)+\sum_{j\not\in {\cJ}^h_{M}}(\lambda^h_{j+1}-
\lambda^h_j)+ (\lambda^h_1-h\alpha)\\ \leq C_1h^{4/3}+KC_1h^{4/3}
+Ch^{-n}h^M + C_1h^{4/3}\,.
\end{multline*}
Taking $M>n+1$, we come to a contradiction.
\end{proof}

\section{Exponential localization of the spectrum}\label{s:weight}
This section is devoted to the proof of Theorem~\ref{t:D}.
Throughout in this section, we will assume that (\ref{e:tr1}) is
satisfied.

\subsection{Weighted $L^2$ spaces}~\\
Let $W$ be an open domain (with regular boundary) in $ M$. Let
\[
b_0(W)=\min \{{\Tr}^+ (B(x))\, :\, x\in W\}.
\]
Denote by $C^{0,1}(\overline{W},\R)$ the class of uniformly
Lipschitz continuous, real-valued functions on $\overline W$.
Introduce the following class of weights
\[
\cW(\overline{W})=\{\Phi\in C^{0,1}(\overline{W},\R) :
\underset{x\in \overline{W}} {\operatorname{ess-inf}} ({\Tr}^+
B(x)-b_0(W)-|\nabla\Phi(x)|^2)
>0\}.
\]
Examples of functions in the class $\cW(\overline{W})$ are given
by the functions $f(x)=(1-\epsilon)d_W(x,X)$, with an arbitrary
$0<\epsilon \leq 1$ and $X\subset W$, where $d_W(x,y)$ is the
distance associated with the (degenerate) Agmon metric
\[
[{\Tr}^+ (B(x))-b_0(W)]_+\cdot { g},
\]
and, for any $x\in \R$, $x_+=\max(x,0)$.

For any $\Phi\in \cW(\overline{W})$ and $h>0$ define the Hilbert
space
\[
L^2_{\Phi/\sqrt{h}}(W)=\{u\in L^2_{loc}(W) : e^{\Phi/\sqrt{h}}u\in
L^2(W)\}
\]
with the norm
\[
\|u\|_{\Phi/\sqrt{h}}=\|e^{\Phi/\sqrt{h}} u\|, \quad u\in
L^2_{\Phi/\sqrt{h}}(W),
\]
where $\|\cdot\|$ denotes the norm in $L^2(W)$:
\[
\|u\|=\left(\int_W|u(x)|^2\,dx\right)^{1/2}, \quad u\in L^2(W).
\]
By $\|\cdot\|_{\Phi/\sqrt{h}}$ we will also denote the norm of a
bounded operator in $L^2_{\Phi/\sqrt{h}}(W)$.

Recall the following important identity (cf. for instance
\cite{HM}).

\begin{lem}\label{l:1.1}
Let $W\subset M$ be an open domain (with $C^2$ boundary) and
$\Phi\in C^{0,1}(\overline{W},\R)$. For any $h>0$, $z\in\C$ and
$u\in \Dom (H^h_W)$ one has
\begin{multline}\label{e:energy}
{\rm Re}\, \int_{W} e^{2\Phi/\sqrt{h}} (H^h_W-z)u \bar{u}\,dx  =
q_W^h(e^{\Phi/\sqrt{h}}u) \\ - h\int_{W} e^{2\Phi/\sqrt{h}}
|\nabla\Phi|^2 |u|^2\,dx - {\rm Re}\, z\int_{W} e^{2\Phi/\sqrt{h}}
|u|^2\,dx.
\end{multline}
\end{lem}

\subsection{Estimates away from the wells}~\\
 Let $W\subset  M$ be a $\Gamma$-invariant
open domain (with a regular boundary). We will start with a slight
extension of \cite[Theorem 3.1]{HM}.

\begin{thm}\label{t:3.1}
There exist constants $C_0>0$ and $h_0>0$ such that for any $h\in
(0,h_0]$ and for any $u\in \Dom (q^h_W)$
\[
h\int_{W} [{\Tr}^+ B(x)-h^{1/4}C_0]\, |u(x)|^2\,dx \leq
(1+h^{1/4}C_0)\, q^h_W(u).
\]
\end{thm}

As a consequence of this theorem, we get
\[
\sigma(H^h_W)\subset [hb_0(W)-C h^{5/4}, +\infty), \quad h\in
(0,h_0],
\]
with some $C>0$ and $h_0>0$.

\begin{prop}\label{p:1}
Let $\Phi\in \cW(\overline{W})$. Assume that $K(h)$ is a bounded
subset in $\C$ such that $K(h)\subset \{z\in \C : {\rm Re}\,z <
h(b_0(W)-\alpha)\}$ for some $\alpha>0$. If $h>0$ is small enough,
then $K(h)\cap \sigma(H_W^h)=\emptyset\,$, and for any $z\in K(h)$
the operator $(H_W^h-z)^{-1}$ defines a bounded operator in
$L^2_{\Phi/\sqrt{h}}(W)$ with
\[
\|(H^h_W-z)^{-1}\|_{\Phi/{\sqrt{h}}}\leq \frac{C}{h}
\]
uniformly on $z\in K(h)$.
\end{prop}

\begin{proof}
By Theorem~\ref{t:3.1} and Lemma~\ref{l:1.1}, for any $z\in \C$,
we have
\begin{multline*}
{\rm Re}\, \int_{W} e^{{2\Phi/\sqrt{h}}} (H^h_W-z)u \bar{u}\,dx
\\ \geq h\int_{W} (1+h^{1/4}C_0)^{-1} [{\Tr}^+
B(x)-h^{1/4}C_0]\,e^{{2\Phi/\sqrt{h}}}\, |u(x)|^2\,dx \\ -
h\int_{W} e^{{2\Phi/\sqrt{h}}}|\nabla\Phi(x)|^2 |u(x)|^2\,dx -
{\rm Re}\,z \int_{W} e^{{2\Phi/\sqrt{h}}} |u(x)|^2\,dx,
\end{multline*}
that implies
\begin{multline*}
{\rm Re}\, \int_{W} e^{{2\Phi/\sqrt{h}}} (H^h_W-z)u \bar{u}\,dx
\\ \geq h\int_{W} e^{{2\Phi/\sqrt{h}}}\, [{\Tr}^+
B(x)-|\nabla\Phi(x)|^2- \frac{{\rm Re}\,z}{h}]\, |u(x)|^2\,dx + c
h^{5/4}\|u\|^2_{\Phi/{\sqrt{h}}} \\ \geq (\alpha +c
h^{1/4})h\|u\|^2_{\Phi/{\sqrt{h}}},
\end{multline*}
and immediately completes the proof.
\end{proof}

\begin{cor}\label{c:h1}
Under the assumptions of Proposition~\ref{p:1}, we have
\begin{multline*}
q_W[e^{{\Phi/\sqrt{h}}}(H^h_W-z)^{-1}v]+
h\|(H^h_W-z)^{-1}v\|^2_{\Phi/{\sqrt{h}}}\\ \leq
\frac{C}{h}\|v\|^2_{\Phi/{\sqrt{h}}},\quad v\in
L^2_{\Phi/\sqrt{h}}(W).
\end{multline*}
\end{cor}

\begin{proof}
By~(\ref{e:energy}), for any $h$ small enough one has
\begin{align*}
q^h[e^{\Phi/\sqrt{h}}(H^h_W-z)^{-1}v] = &{\rm Re}\,(
e^{2\Phi/\sqrt{h}} v, (H^h_W-z)^{-1}v)\\ & + h\| |\nabla\Phi|
(H^h_W-z)^{-1}v\|_{\Phi/{\sqrt{h}}}^2\\ & + {\rm Re}\, z
\|(H^h_W-z)^{-1}v\|_{\Phi/{\sqrt{h}}}^2.
\end{align*}
Now we know that ${\rm Re}\,z < h(b_0(W)-\alpha)$, $|\nabla\Phi|$
is uniformly bounded and
\begin{align*}
{\rm Re}\,( e^{2\Phi/\sqrt{h}} v, (H^h_W-z)^{-1}v) & \leq
\frac{1}{2}\left(\frac{1}{h}\|v\|^2_{\Phi/{\sqrt{h}}} +
h\|(H^h_W-z)^{-1}v\|_{\Phi/{\sqrt{h}}}^2\right)\\ & \leq
\frac{C}{h}\|v\|^2_{\Phi/{\sqrt{h}}},
\end{align*}
that completes the proof.
\end{proof}

\subsection{Estimates near the wells}~\\
In this section, we will assume that $W$ is a relatively compact
domain (with smooth boundary) in $\cF$ such that
\[
U_{\epsilon_1}=\{x\in \cF\,:\, {\Tr}^+ (B(x)) < b_0+\epsilon_1\}
\]
is contained in $W$ for some $\epsilon_1<\epsilon_0$.

\begin{prop}\label{p:10}
Assume that $K(h)$ is a bounded subset in $\C$ such that
$K(h)\subset \{z\in \C : {\rm Re}\,z < h(b_0+\epsilon_1)\}$ and,
if $h>0$ is small enough, then, for any $\epsilon >0$,
\[
{\rm dist}\,(K(h), \sigma(H^h_W))\geq \frac{1}{C_\epsilon}
e^{-\epsilon/\sqrt{h}}.
\]
Let $\Phi\in \cW(\overline{W})$ such that $\Phi\equiv 0$ on
$U_{\epsilon_1}$. Then for any $z\in K(h)$ the operator
$(H_W^h-z)^{-1}$ defines a bounded operator in
$L^2_{\Phi/\sqrt{h}}(W)$ and for any $\epsilon>0$
\[
\|(H_W^h-z)^{-1}\|_{\Phi/{\sqrt{h}}}\leq
C_{1,\epsilon}e^{\epsilon/\sqrt{h}}\,.
\]
\end{prop}

\begin{proof}
For every sufficiently small $\eta>0$, take any $\chi_{1,\eta}\in
C^\infty_c(W)$ such that $\chi_{1,\eta}\equiv 1$ in a neighborhood
of $\{x\in W:\Phi(x)\leq 2\eta \}$, $\Phi\leq 3\eta$ on $\supp
\chi_{1,\eta}$. Let $\chi^\prime_{1,\eta}\in C^\infty(W)$,
$\chi^\prime_{1,\eta}\geq 0$ satisfy $(\chi_{1,\eta})^2 +
(\chi^\prime_{1,\eta})^2=1$. We can assume that there exists a
constant $C$ such that, for all sufficiently small $\eta>0$,
\[
\eta (|\nabla \chi_{1,\eta}|+|\nabla \chi^\prime_{1,\eta}|)\leq C\,.
\]
Then we have
\begin{multline*}
q_W(e^{{\Phi/\sqrt{h}}}u)= q(\chi_{1,\eta} e^{{\Phi/\sqrt{h}}}u)
+q(\chi^\prime_{1,\eta} e^{{\Phi/\sqrt{h}}}u)\\ -h^2\| |\nabla
\chi_{1,\eta}| e^{{\Phi/\sqrt{h}}} u\|^2-h^2\| |\nabla
\chi^\prime_{1,\eta}| e^{{\Phi/\sqrt{h}}} u\|^2\,.
\end{multline*}
By (\ref{e:energy}), it follows that
\begin{multline}
\label{e:energy1} q_W(\chi^\prime_{1,\eta} e^{{\Phi/\sqrt{h}}}u)\\
- h\int_{W} e^{{2\Phi/\sqrt{h}}} |\nabla\Phi|^2
|\chi^\prime_{1,\eta}u|^2\,dx - {\rm Re}\, z\int_{W}
e^{{2\Phi/\sqrt{h}}} |\chi^\prime_{1,\eta} u|^2\,dx\\ - h^2\|
|\nabla \chi_{1,\eta}| e^{{\Phi/\sqrt{h}}}\chi^\prime_{1,\eta}
u\|^2  - h^2\| |\nabla \chi^\prime_{1,\eta}|
e^{{\Phi/\sqrt{h}}}\chi^\prime_{1,\eta} u\|^2
\\ =  {\rm Re}\, \int_{W} e^{{2\Phi/\sqrt{h}}} (H^h_W-z)u
\bar{u}\,dx - q_W(\chi_{1,\eta} e^{{\Phi/\sqrt{h}}}u) \\  +
h\int_{W} e^{{2\Phi/\sqrt{h}}} |\nabla\Phi|^2
|\chi_{1,\eta}u|^2\,dx  + {\rm Re}\, z\int_{W}
e^{{2\Phi/\sqrt{h}}} |\chi_{1,\eta} u|^2\,dx\\ + h^2\| |\nabla
\chi_{1,\eta}| e^{{\Phi/\sqrt{h}}}\chi_{1,\eta} u\|^2+ h^2\|
|\nabla \chi^\prime_{1,\eta}| e^{{\Phi/\sqrt{h}}}\chi_{1,\eta}
u\|^2\,.
\end{multline}
Put $\eta=\alpha \sqrt{h}$ with sufficiently large $\alpha>0$.
Taking into account the fact that
\[
|\nabla \chi_{1,\eta}|+|\nabla \chi^\prime_{1,\eta}|\leq
{C/\eta}={C/\alpha \sqrt{h}},
\]
we get the following estimate for
the right-hand side of (\ref{e:energy1})
\begin{multline*}
h\int_{W} e^{{2\Phi/\sqrt{h}}} |\nabla\Phi|^2
|\chi_{1,\eta}u|^2\,dx + {\rm Re}\, z\int_{W} e^{{2\Phi/\sqrt{h}}}
|\chi_{1,\eta} u|^2\,dx \\ + h^2\| |\nabla \chi_{1,\eta}|
e^{{\Phi/\sqrt{h}}}\chi_{1,\eta} u\|^2 + h^2\| |\nabla
\chi^\prime_{1,\eta}| e^{{\Phi/\sqrt{h}}}\chi_{1,\eta} u\|^2
\\ \leq C h \|e^{{\Phi/\sqrt{h}}} \chi_{1,\eta}u\|^2\,.
\end{multline*}
From the other side, proceeding as in the proof of
Proposition~\ref{p:1} and using Theorem~\ref{t:3.1}, we get the
estimate for the left-hand side of (\ref{e:energy1})
\begin{multline*}
q_W(\chi^\prime_{1,\eta} e^{{\Phi/\sqrt{h}}}u) \\ - h\int_{W}
e^{{2\Phi/\sqrt{h}}} |\nabla\Phi|^2 |\chi^\prime_{1,\eta}u|^2\,dx
- {\rm Re}\, z\int_{W} e^{{2\Phi/\sqrt{h}}} |\chi^\prime_{1,\eta}
u|^2\,dx \\ - h^2\| |\nabla \chi_{1,\eta}|
e^{{\Phi/\sqrt{h}}}\chi^\prime_{1,\eta} u\|^2 - h^2\| |\nabla
\chi^\prime_{1,\eta}| e^{{\Phi/\sqrt{h}}}\chi^\prime_{1,\eta}
u\|^2 \\ \geq h\int_{W} e^{{2\Phi/\sqrt{h}}}\, \left[{\Tr}^+
B(x)-|\nabla\Phi(x)|^2- \frac{{\rm
Re}\,z}{h}-\frac{1}{\alpha^2}\right]\, |\chi^\prime_{1,\eta}
u(x)|^2\,dx\\ + c h^{5/4}\|e^{{\Phi/\sqrt{h}}}
\chi^\prime_{1,\eta}u\|^2 \geq C h \|e^{{\Phi/\sqrt{h}}}
\chi^\prime_{1,\eta}u\|^2\,.
\end{multline*}
Thus we get the estimate
\[
c h \|e^{{\Phi/\sqrt{h}}} u\|^2 \leq {\rm Re}\, \int_{W}
e^{{2\Phi/\sqrt{h}}}(H^h_W-z)u \bar{u}\,dx \\ + C h
\|e^{{\Phi/\sqrt{h}}} \chi_{1,\eta}u\|^2\,.
\]
It remains to show that, for any $\epsilon>0$
\begin{equation*}
\|e^{{\Phi/\sqrt{h}}} \chi_{1,\eta}u\| \leq C_\epsilon
e^{{\epsilon/\sqrt{h}}} \|e^{{\Phi/\sqrt{h}}}(H^h_W-z)u\|\,,
\end{equation*}
or equivalently,
\begin{equation}\label{e:chi}
\|\chi_{1,\eta}(H^h_W-z)^{-1}u\|_{\Phi/{\sqrt{h}}} \leq C_\epsilon
e^{{\epsilon/\sqrt{h}}} \|u\|_{\Phi/{\sqrt{h}}}, \quad u\in
L^2_{\Phi/{\sqrt{h}}}(W)\,.
\end{equation}
For this, we choose a function $\chi_{2,\eta}\in C^\infty_c(W)$
such that $\chi_{2,\eta}\equiv 1$ in a neighborhood of $\{x\in
W:\Phi(x)\leq \eta\}$, $\Phi\leq 2\eta$ on $\supp \chi_{2,\eta}$.
In particular, $\chi_{1,\eta}\equiv 1$ on $\supp \chi_{2,\eta}$.
We can assume that there exists a constant $C$ such that for all
sufficiently small $\eta>0$
\begin{equation}\label{e:chi2}
\eta |\nabla \chi_{2,\eta}|+\eta^2|\Delta \chi_{2,\eta}|\leq C.
\end{equation}

Let $M_0=\{x\in W: \Phi(x)\geq 2\eta\}$. Then we have
\begin{multline*}
(H^h_{W}-z)^{-1}u=(1-\chi_{2,\eta})(H^h_{M_0}-z)^{-1}(1-\chi_{1,\eta})u+
(H^h_{W}-z)^{-1}\chi_{1,\eta}u \\ + (H^h_{W}-z)^{-1}\chi_{1,\eta}
[H^h_W,\chi_{2,\eta}] (H^h_{M_0}-z)^{-1}(1-\chi_{1,\eta})u\,.
\end{multline*}
We consider three terms in the right hand side of the last
identity separately. For the first one we use
Proposition~\ref{p:1} and obtain
\begin{equation}\label{e:1}
\|\chi_{1,\eta}(1-\chi_{2,\eta})(H^h_{M_0}-z)^{-1}
(1-\chi_{1,\eta})u\|_{\Phi/{\sqrt{h}}}\leq \frac{C}{h}
\|u\|_{\Phi/{\sqrt{h}}}\,.
\end{equation}
For the second term, since $\Phi\leq 3\eta$ on $\supp
\chi_{1,\eta}$, we have
\[
\|\chi_{1,\eta}(H^h_{W}-z)^{-1}\chi_{1,\eta}u\|_{\Phi/{\sqrt{h}}}
\leq e^{3\alpha} \|(H^h_{W}-z)^{-1}\chi_{1,\eta} u\|\,.
\]
By the assumptions and the fact that $\Phi\geq 0$, it follows that
\[
\|(H^h_{W}-z)^{-1}\chi_{1,\eta} u\|\leq e^{\epsilon/\sqrt{h}}
\|\chi_{1,\eta} u\| \leq C_1 e^{{\epsilon/\sqrt{h}}} \|
u\|_{\Phi/{\sqrt{h}}}\,.
\]
So we get for the second term
\begin{equation}\label{e:2}
\|\chi_{1,\eta}(H^h_{W}-z)^{-1}\chi_{1,\eta}u\|_{\Phi/{\sqrt{h}}}
\leq C_2 e^{{\epsilon/\sqrt{h}}} \| u\|_{\Phi/{\sqrt{h}}}\,.
\end{equation}

For the third term we put $w=(H^h_{M_0}-z)^{-1} (1-
\chi_{1,\eta})u$. By (\ref{e:2}), it follows that
\[
\|\chi_{1,\eta}(H^h_{W}-z)^{-1}\chi_{1,\eta} [H^h_W,
\chi_{2,\eta}] w\|_{\Phi/{\sqrt{h}}}\leq C_1
e^{{\epsilon/\sqrt{h}}}\|[H^h_W, \chi_{2,\eta}]
w\|_{\Phi/{\sqrt{h}}}\,.
\]
Now we have
\[
[H^h_W, \chi_{2,\eta}]w=2i h\, d \chi_{2,\eta}\cdot (ih\,d+{\bf
A})w + h^2\Delta \chi_{2,\eta} w\,.
\]
Therefore, taking into account (\ref{e:chi2}), we get
\begin{align*}
\|[H^h_W, \chi_{2,\eta}] w\|^2_{\Phi/{\sqrt{h}}} & \leq C
(h\|(ih\,d+{\bf A})w\|^2_{\Phi/{\sqrt{h}}}+ h^2
\|w\|^2_{\Phi/{\sqrt{h}}})\\ & \leq C (h q_W[e^{{\Phi/\sqrt{h}}}w]
+ h^2 \|w\|^2_{\Phi/{\sqrt{h}}})\,.
\end{align*}
By Corollary~\ref{c:h1}, we have
\begin{align*}
\|[H^h_W, \chi_{2,\eta}] w\|^2_{\Phi/{\sqrt{h}}}  \leq & C (h
q_W[e^{{\Phi/\sqrt{h}}}(H^h_{M_0}-z)^{-1} (1- \chi_{1,\eta})u]\\ &
+ h^2 \|(H^h_{M_0}-z)^{-1} (1-
\chi_{1,\eta})u\|^2_{\Phi/{\sqrt{h}}}) \\ \leq & C \|(1-
\chi_{1,\eta})u\|^2_{\Phi/{\sqrt{h}}} \leq C
\|u\|^2_{\Phi/{\sqrt{h}}}\,.
\end{align*}
So we get for the third term
\begin{multline}\label{e:3}
\|\chi_{1,\eta}(H^h_{W}-z)^{-1}\chi_{1,\eta} [H^h_W,
\chi_{2,\eta}] (H^h_{M_0}-z)^{-1} (1-
\chi_{1,\eta})u\|_{\Phi/{\sqrt{h}}}\\ \leq C_{3,\epsilon}
e^{{\epsilon/\sqrt{h}}} \|u\|_{\Phi/{\sqrt{h}}}\,.
\end{multline}
Now (\ref{e:chi}) follows by adding the estimates (\ref{e:1}),
(\ref{e:2}) and (\ref{e:3}).
\end{proof}

\begin{cor}\label{c:h2}
Under the assumptions of Proposition~\ref{p:10}, we have, for any
$\epsilon>0$\,,
\begin{multline*}
q_W[e^{{\Phi/\sqrt{h}}}(H^h_W-z)^{-1}v]+
h\|(H^h_W-z)^{-1}v\|^2_{\Phi/{\sqrt{h}}}\\ \leq
C_{2,\epsilon}e^{\epsilon/\sqrt{h}}
\|v\|^2_{\Phi/{\sqrt{h}}},\quad v\in L^2_{\Phi/\sqrt{h}}(W)\,.
\end{multline*}
\end{cor}

\subsection{Proof of Theorem~\ref{t:D}}~\\
Let us assume that (\ref{e:tr1}) and (\ref{e:trB}) are satisfied.
We have
\[
\{x\in \cF\,:\, {\Tr}^+ (B(x)) < b_0+\epsilon_2\} = U_{\epsilon_1}
= \bigcup_{j=1}^N U_{j,\epsilon_2}\,,
\]
where $U_{j,\epsilon_2}\subset \cF, j=1,2,\cdots, N,$ are
relatively compact, connected and pairwise disjoint domains such
that $U_{j,\epsilon_2}\cap \partial\cF=\emptyset$\,. Let
$M_j=\overline{U_{j,\epsilon_2}}$\,,  $j=1,2,\cdots, N$.
Theorem~\ref{t:D} follows immediately from the following

\begin{prop}
Assume that $K(h)$ is a bounded subset in $\C$ such that
$K(h)\subset \{z\in \C : {\rm Re}\,z < h(b_0+\epsilon_1)\}$ and,
if $h>0$ is small enough, then, for any $\epsilon >0$\,,
\[
{\rm dist}\,(K(h), \sigma(H^h_{M_j}))\geq
\frac{1}{C_\epsilon}e^{-\epsilon/\sqrt{h}}, \quad j=1,2,\cdots, N\,.
\]
Then, for any $h>0$ small enough, $K(h)\cap \sigma(H^h) =
\emptyset$\,.
\end{prop}

\begin{proof}
Take any $\eta>0$ such that $\epsilon_1+3\eta<\epsilon_2$. Let
\[
M_0=M\setminus \bigcup_{\gamma\in\Gamma}\bigcup_{j=1}^N \gamma (
U_{j,\epsilon_1+\eta})= \{x\in { M}\,:\, {\Tr}^+ (B(x)) \geq
b_0+\epsilon_1+\eta \}\,.
\]
Take any function $\phi_j\in C^\infty_c({ M})$ such that $\supp
\phi_j \subset U_{j,\epsilon_1+2\eta}$, $\phi_j\equiv 1$ on
$U_{j,\epsilon_1+\eta}$. Let
\[
\phi_0=1-\sum_{\gamma\in\Gamma}\sum_{j=1}^N\gamma^*\phi_j\,.
\]
Then $\supp \phi_0\subset M_0$. Let $\psi_j\in C^\infty_c({ M})$,
$(j=1,2,\cdots, N$, such that $\supp \psi_j \subset
U_{j,\epsilon_1+3\eta}$, $\psi_j\equiv 1$ in a neighborhood of
$U_{j,\epsilon_1+2\eta}$. Take any $\Gamma$-periodic function
$\psi_0\in C^\infty({ M})$ such that $\supp \psi_0 \subset M_0$,
$\psi_0\equiv 1$ in a neighborhood of $M\setminus
\cup_{\gamma\in\Gamma}\cup_{j=1}^N \gamma
(U_{j,\epsilon_1+2\eta})$. In particular, we have
$\phi_j\psi_j=\phi_j$, for $ j=0,1,2,\cdots, N$.

Recall that the magnetic translations $T_\gamma, \gamma\in\Gamma,$
are unitary operators in $L^2({ M})$, which commute with the
periodic magnetic Schr\"odinger operator $H^h$:
\[
T_\gamma H^h=H^hT_\gamma, \quad \gamma\in\Gamma,
\]
and each $T_\gamma$ takes $L^2(\cF)$ to $L^2(\gamma\cF)$ (see for
instance \cite{MS,KMS} and references therein for more details).
They satisfy
$$ T_e={\id}, \quad T_{\gamma_1} T_{\gamma_2} =
{\sigma}(\gamma_1,\gamma_2) T_{\gamma_1 \gamma_2}, \quad \gamma_1,
\gamma_2 \in \Gamma\,.$$
Here $\sigma$ is a $2$-cocycle on $\Gamma$
i.e. $\sigma:\Gamma\times\Gamma\to U(1)$ such that
\begin{align*}
{\sigma}(\gamma,e)  = & {\sigma}(e,\gamma)=1,\quad
\gamma\in\Gamma;\\ {\sigma}(\gamma_1,\gamma_2)
{\sigma}(\gamma_1\gamma_2, \gamma_3) = &
{\sigma}(\gamma_1,\gamma_2\gamma_3) {\sigma}(\gamma_2,\gamma_3),
\quad \gamma_1, \gamma_2, \gamma_3\in \Gamma\,.
\end{align*}

For any $h>0$ small enough and any $z\in K(h)$, define a bounded
operator $R^h(z)$ in $L^2( M)$ as
\[
R^h(z)=\sum_{j=1}^N\sum_{\gamma\in\Gamma} T_\gamma \psi_j
(H^h_{M_j}-z)^{-1}\phi_j T^*_\gamma + \psi_0
(H^h_{M_0}-z)^{-1}\phi_0\,.
\]
Then
\[
(H^h-z)R^h(z)=I-K^h(z)\,,
\]
where
\[
K^h(z)=\sum_{j=1}^N\sum_{\gamma\in\Gamma} T_\gamma [H^h,\psi_j]
(H^h_{M_j}-z)^{-1}\phi_j T^*_\gamma + [H^h,\psi_0]
(H^h_{M_0}-z)^{-1}\phi_0\,.
\]

\begin{lem}\label{l:kh}
There exist $C,c >0$ such that, for any $h>0$ small enough and
$z\in K(h)$,  the operator $K^h(z)$ defines a bounded operator in
$L^2( M)$ with the norm estimate
\[
\|K^h(z)\|\leq C e^{-c/\sqrt{h}}\,.
\]
\end{lem}

\begin{proof}
For any $j=1,2,\cdots, N$, consider a weight function $\Phi_j\in
\cW(M_j)$ given by
$\Phi_j(x)=d_{U_{j,\epsilon_1+2\eta}}(x,U_{j,\epsilon_1+2\eta})$.
By construction, $\Phi_j(x)\geq c_j>0$ on $\supp d\psi_j$,
$\Phi_j(x)\equiv 0$ on $\supp \phi_j$. For any $w\in {\rm
Dom}\,H^h$, we have
\[
[H^h, \psi_j]w=2i h\, d \psi_j\cdot (ih\,d+{\bf A})w + h^2\Delta
\psi_j\, w\,.
\]
This  implies the estimate
\begin{align*}
\|[H^h, \psi_j] w\|^2_{\Phi_j/{\sqrt{h}}} & \leq C (h\|(ih\,d+{\bf
A})w\|^2_{\Phi_j/{\sqrt{h}}}+ h^2 \|w\|^2_{\Phi_j/{\sqrt{h}}})\\ &
\leq C (h q_{M_j}[e^{{\Phi_j/\sqrt{h}}}w] + h^2
\|w\|^2_{\Phi_j/{\sqrt{h}}})\,.
\end{align*}
Therefore, for any $u\in L^2( M)$, we obtain
\begin{multline*}
\|[H^h,\psi_j] (H^h_{M_j}-z)^{-1}\phi_j u\|^2_{L^2( M)}\\
\begin{aligned}
= & \|[H^h,\psi_j] (H^h_{M_j}-z)^{-1}\phi_j u\|^2_{L^2(M_j)}\\
\leq & e^{-c_j/\sqrt{h}} \|[H^h,\psi_j] (H^h_{M_j}-z)^{-1}\phi_j
u\|^2_{\Phi_j/\sqrt{h}}\\  \leq & C e^{-c_j/\sqrt{h}} (h
q_{M_j}[e^{{\Phi_j/\sqrt{h}}}(H^h_{M_j}-z)^{-1}\phi_j u]\\ & + h^2
\|(H^h_{M_j}-z)^{-1}\phi_j u\|^2_{\Phi_j/{\sqrt{h}}})\,.
\end{aligned}
\end{multline*}
It follows from Corollary \ref{c:h2} that, for any $\epsilon>0$,
\begin{align*}
\|[H^h,\psi_j] (H^h_{M_j}-z)^{-1}\phi_j u\|_{L^2( M)} \leq &
C_\epsilon e^{-(c_j-\epsilon)/\sqrt{h}} \|\phi_j
u\|_{\Phi_j/{\sqrt{h}}}\\ = & C_\epsilon
e^{-(c_j-\epsilon)/\sqrt{h}} \|\phi_j u\|_{L^2(M_j)}\\ = &
C_\epsilon e^{-(c_j-\epsilon)/\sqrt{h}} \| u\|_{L^2( M)}\,.
\end{align*}
Similarly, using Corollary \ref{c:h1}, one can get
\[
\|[H^h,\psi_0] (H^h_{M_0}-z)^{-1}\phi_0 u\|_{L^2( M)} \leq C_0
e^{-c_0/\sqrt{h}} \| u\|_{L^2( M)}\,.
\]
Taking into account that the sets $\gamma(\supp \phi_j)$ with
$j=1,2,\cdots,N$ and $\gamma\in \Gamma$ are disjoint, we get
\begin{align*}
\|K^h(z)u\| & \leq C e^{-c/\sqrt{h}} (\sum_{j=1}^N
\sum_{\gamma\in\Gamma} \|\phi_j T^*_\gamma u\| + \|\phi_0 u\|)\\ &
\leq C_1 e^{-c/\sqrt{h}} \|u\|\,.
\end{align*}
This  completes the proof.
\end{proof}

It follows from Lemma~\ref{l:kh} that, for all sufficiently small
$h>0$ and $z\in K(h)$, the operator $I+K^h(z)$ is invertible in
$L^2( M)$. Then the operator $H^h-z$ is invertible in $L^2( M)$
with
\[
(H^h-z)^{-1}= R^h(z)(I-K^h(z))^{-1}\,,
\]
and $K(h)\cap \sigma(H^h)=\emptyset$ as desired.
\end{proof}

\end{document}